\documentclass[12pt]{amsart}
\usepackage{amssymb,amsmath,amsfonts}
\usepackage{enumerate}
\usepackage[letterpaper, margin=1.2in]{geometry}

\allowdisplaybreaks

\newcommand{\N}{\mathbb N}
\newcommand{\R}{\mathbb R}

\newcommand{\rar}{\rightarrow}
\newcommand{\mc}{\mathcal}

\DeclareMathOperator{\card}{card}

\newtheorem*{theorem*}{Theorem}

\begin{document}

\title[A measure-theoretic result for approximation by Delone sets]{A
  measure-theoretic result for\\[3mm] approximation by Delone sets}
\author{Michael Baake, Alan Haynes}


\begin{abstract}
  With a view to establishing measure-theoretic approximation
  properties of Delone sets, we study a setup which arises naturally
  in the problem of averaging almost periodic functions along
  exponential sequences. In this setting, we establish a full converse
  of the Borel--Cantelli lemma. This provides an analogue of more
  classical problems in the metric theory of Diophantine
  approximation, but with the distance to the nearest integer function
  replaced by distance to an arbitrary Delone set.
\end{abstract}

\maketitle

\section{Introduction}
In a recent work \cite{BaakGrim2017}, motivated by problems emerging
in aperiodic order \cite{BaakGrim2013}, we considered the problem of
establishing asymptotic formulas for averages of Bohr almost periodic
functions along exponential sequences of the form
$(\alpha^nx)_{n\in\N}$, where $\alpha$ is a fixed real number with
$|\alpha|>1$, and $x\in\R$ is arbitrary. The almost everywhere
convergence results which were obtained can be viewed as analogues of
Birkhoff's ergodic theorem \cite{Walt1982}. Their proofs relied on
knowing that, for almost every $x\in\R$, the elements of our
exponential sequence with $n\le N$ cannot be approximated too well by
elements of a given Delone set. In fact, it was precisely this
occurrence of a Delone set that showed how natural the replacement of
the integers by such a more general set is, and how well it blends
into the realm of metric Diophantine approximation problems. Once
again, systems of aperiodic order have thus suggested how one can go
meaningfully beyond the standard lattice setting, while remaining
firmly within the realm of physical and mathematical significance. In
this short note, we return to this problem in order to more fully
investigate the Diophantine approximation properties of Delone
sets. Our main result is a generalization of Khintchine's theorem
(even without monotonicity) to this setting (see \cite{Buge2012} or
\cite{Harm1998} for detailed background on the metric theory of
Diophantine approximation). It provides a characterization which tells
us precisely when we should expect, for almost every $x\in\R$, to have
infinitely many approximations of a certain quality.\vspace*{.1in}

Recall that a Delone set is a subset $Y$ of a metric space $(X,d)$ for
which there exist positive constants $r$ and $R$ such that:
\begin{enumerate}\itemsep=2pt
\item[(i)] For any $y,y'\in Y$ with $y\not=y'$, we have that
  $d(y,y')\ge 2r$, and \vspace*{2pt}
\item[(ii)] For any $x\in X$, there exists a $y\in Y$ such that
  $d(x,y)\le R$.
\end{enumerate}
The supremum over all $r$ satisfying (i) is called the packing radius
of $Y$, and the infimum over all $R$ satisfying (ii) is called the
covering radius. See \cite[Section 2.1]{BaakGrim2013} for more
background on the basic theory of Delone sets.\vspace*{.1in}

In what follows, we suppose that $Y$ is a Delone set in $\R$ with
packing radius $r$ and covering radius $R$. For $x\in\R$ and $\rho>0$,
we use $B(x,\rho)$ to denote the closed ball of radius $\rho$ centered
at $x$. The symbol $\lambda$ denotes Lebesgue measure on $\R$, and
$\card (S)$ denotes the cardinality of a set $S$. We also use the
standard Vinogradov $\ll$ notation, so that if
$f,g:D\longrightarrow\R$ are two functions on some common domain $D$,
then $f\ll g$ means that there exists a constant $C>0$ such that
$|f(x)|\le C|g(x)|$ for all $x\in D$. Our main result is the following
theorem.

\begin{theorem*}
  Let\/ $\psi:\N\longrightarrow [0,\infty)$ be any function, and
  suppose that\/ $\alpha$ is a real number with\/ $|\alpha|>1$. If
\begin{equation}\label{eqn.PsiConv}
\sum_{n=1}^\infty\psi (n)<\infty
\end{equation}
then, for almost every\/ $x\in\R$, there are only finitely many\/
$n\in\N$ and\/ $y\in Y$ for which
\begin{equation}\label{eqn.PsiIneq}
|\alpha^nx-y|\le \psi (n).
\end{equation}
On the other hand, if
\begin{equation}\label{eqn.PsiDiv}
\sum_{n=1}^\infty\psi (n)=\infty ,
\end{equation}
then, for almost every\/ $x\in\R$, there are infinitely many\/
$n\in\N$ and\/ $y\in Y$ that satisfy the above inequality.
\end{theorem*}

We point out that our result amounts to a full converse of the
Borel--Cantelli lemma for special collections of subsets of $\R$. The
less trivial direction of our proof is the divergence part, in which
we assume that \eqref{eqn.PsiDiv} holds, and prove that the associated
limsup set (the set of real numbers $x$ for which the relevant
inequalities have infinitely many solutions) has full measure. It is
quite common in these types of problems to first establish a zero-full
lemma (e.g.\ \cite{BereDickVela2006,Gall1961, Harm1998}), which proves
that, whether or not \eqref{eqn.PsiDiv} holds, the associated limsup
set either has zero measure, or its complement does. Then all that is
left to show, for the difficult part of the proof, is that the limsup
set has {\em positive} measure. In our proof, we deviate slightly from
this approach and, instead of establishing a separate zero-full lemma,
we use a simple application of the Lebesgue density theorem
\cite[Theorem 3.21]{Foll1984} to complete our proof. In essence, what
we show is that, at all scales, the intersections of the individual
intervals which contribute to the limsup set behave in the same
way. In other words, roughly the same picture appears after zooming in
on any part of the real line. Since the integers in our setup have
been replaced by a completely arbitrary Delone set, this is clearly
another manifestation of a behavior which should be subsumed under the
theme of `aperiodic order'.

\section{Proof of main theorem}
We will consider only the case when $\alpha>1$, since the case when
$\alpha<-1$ follows as an easy corollary. Let $a<b$ be real numbers
and assume, with little loss of generality, that $a\ge 0$ (the cases
with $a<0$ can be dealt with by trivial modifications of our
argument). For each $n\in\N$, define a set $\mc{A}_n\subseteq [a,b)$
by
\[\mc{A}_n=\{x\in [a,b): |\alpha^nx-y|\le \psi(n)~\text{for some}~y\in Y\}.\]
In what follows, all constants implied by the use of the $\ll$
notation will be universal, not depending on $a,b,r,$ or $R$, unless
otherwise stated. Whenever implied constants do depend on some of
these quantities, we will indicate this by attaching the appropriate
subscripts to the $\ll$ symbol.\vspace*{.1in}

The proof of the first part of the theorem is a straightforward
application of the convergence part of the Borel--Cantelli
lemma. Write the non-negative elements of $Y$ in increasing order as
\[y_1<y_2<\cdots\]
and, for each $n\in\N$, define
\[\mc{I}_n=\left\{i\in\N:[a,b)\cap B\left(\frac{y_i}{\alpha^n},
      \frac{\psi(n)}{\alpha^n}\right)\not=\varnothing\right\}.\]
From our hypothesis on $Y$ we have, for all sufficiently large $n$
(depending on $a$ and $b$), that
\[\frac{(b-a)\,\alpha^n}{R}\ll\card(\mc{I}_n)\ll\frac{(b-a)\,\alpha^n}{r},\]
and it follows from this that, for $n$ sufficiently large,
\[\lambda(\mc{A}_n)\ll (b-a)\,\psi (n)/r.\]
Since we are assuming that \eqref{eqn.PsiConv} holds, we conclude that
almost every $x\in [a,b)$ falls in only finitely many of the sets
$\mc{A}_n$, which is equivalent to the assertion that there are only
finitely many solutions to the inequality \eqref{eqn.PsiIneq}. This
gives the conclusion of the first part of the theorem.\vspace*{.1in}

The proof of the second part is slightly more complicated. Ideally, we
would like to demonstrate that the sets $\mc{A}_n$ above are
quasi-independent, or in other words that
\[
    \lambda(\mc{A}_m\cap\mc{A}_n)\ll_{r,R} (b-a)^{-1}
    \lambda(\mc{A}_m)\lambda(\mc{A}_n)
    \quad\text{for}\quad m\not=n.
\]
Unfortunately, this is not quite true, so we need some technical
modifications in our setup. First of all, choose an integer $J$ with
the property that
\begin{equation}\label{eqn.JChoice}
\alpha^J\ge\frac{2R(1+\tfrac{2}{r})}{r}+1,
\end{equation}
and then choose a residue class $j\in \{0,1,\ldots ,J\! -\! 1\}$
modulo $J$ with the property that
\[\sum_{n=1}^\infty\psi (nJ+j)=\infty.\]
Since there are only $J$ residue classes to choose from, and since we
are assuming that \eqref{eqn.PsiDiv} holds, it is clear that this is
possible.\vspace*{.1in}

Now write $\beta=\alpha^J$ and, for each $n\in\N$, define a set
$Y^{(n)}\subseteq Y$ by
\begin{equation}\label{eqn.Y^nDef}
Y^{(n)}=Y\setminus \bigcup_{m=1}^{n-1}B(\beta^{n-m}Y,1).
\end{equation}
Here we are using $B(\beta^{n-m}Y,1)$ to denote the union over
$y\in Y$ of the balls $B(\beta^{n-m}y,1).$ The reason for introducing
the sets $Y^{(n)}$, which will become clearer later in the proof, is
to remove the bad overlaps which occur between the sets $\mc{A}_n$
from the previous argument. However, we will still need to show that
we have not discarded too much from $Y$ so as to make the sum of the
measures of our new sets fail to diverge. For each $X\in\R$, we have
that
\[\card\{y\in Y:aX\le y< bX\}\ge \left\lfloor\frac{(b-a)X}{R}\right\rfloor\]
and, for each $\ell\in\N$, we also have that
\[\card\{y\in\beta^\ell Y:aX\le y< bX\}\le\frac{(b-a)X}{r\beta^\ell}+1.\]
Since the number of points of $Y$ in a ball of radius $1$ is bounded
above by $1+\tfrac{2}{r},$ it follows that, for $n,N\in\N$,
\begin{align}
  \card\{y\in Y^{(n)}:aX\le y< bX\}
  &  \ge\left\lfloor\frac{(b-a)X}{R}\right
    \rfloor-\sum_{\ell=1}^{n-1}
    \left(\frac{(b-a)X}{r\beta^\ell}+
    1\right)\cdot(1+\tfrac{2}{r}) \nonumber \\
  &\ge (b-a)X\left(\frac{1}{R}-
    \frac{1+\tfrac{2}{r}}{r(\beta-1)}\right)
    -1-n(1+\tfrac{2}{r})\nonumber \\[2mm]
  &\ge \frac{(b-a)X}{2R}-1-n(1+\tfrac{2}{r}).\label{eqn.Y^nLowerBd}
\end{align}
The final inequality here is a result of our choice of $J$ in
\eqref{eqn.JChoice}.\vspace*{.1in}

For $n\in\N$ we now define $\mc{A}_n'\subseteq [0,1)$, our replacement
for the set $\mc{A}_n$ from above, by
\[\mc{A}_n'=\{x\in [a,b): |\alpha^j\beta^nx-y|\le
  \psi(nJ+j)~\text{for some}~y\in Y^{(n)}\}.\] From
\eqref{eqn.Y^nLowerBd}, and using the arguments from above, we have
that, for all sufficiently large $n$ (depending again only on $a$ and
$b$),
\[
    \lambda(\mc{A}_n')\gg \frac{(b-a)\, \psi(nJ+j)}{R}.
\]
For $m\not= n$, we now would like to derive an upper bound for the
measure of the intersection of $\mc{A}_m'$ with $\mc{A}_n'$. With this
purpose in mind,
let
\[\delta=\delta(m,n)=\min\left\{\frac{2\psi(mJ+j)}
    {\alpha^j\beta^m},\frac{2\psi(nJ+j)}{\alpha^j\beta^n}\right\},\]
and
\[
   \Delta=\Delta(m,n)=\max\left\{\frac{2\psi(mJ+j)}{\alpha^j\beta^m},
       \frac{2\psi(nJ+j)}{\alpha^j\beta^n}\right\}.\vspace{1mm}
\]
Each set $\mc{A}_n'$ is a union of connected components, which we
refer to as its component intervals. The component intervals of
$\mc{A}_n'$ are intervals centered at points of the form
$y/\alpha^j\beta^n$. To be fully accurate, at the endpoints of $[a,b)$
it may be the case that there are (at most $2$) component intervals
which are not of this form, but this fact is negligible in the
argument we are about to give. If a component interval from
$\mc{A}_m'$ intersects a component interval from $\mc{A}_n'$, the
centers of these intervals must be within $\Delta$ of one
another. Furthermore, the intersection of any two such intervals has
measure at most $\delta$. This translates into the upper bound
\begin{align}
  \lambda(\mc{A}_m'\cap\mc{A}_n')\ll\delta\cdot
  \card\left\{(y,y')\in Y^{(m)}\!\times Y^{(n)}~:
  \phantom{\frac{1}{1}}\right.
  &a\le \frac{y}{\alpha^j\beta^m}< b,~ a\le
    \frac{y'}{\alpha^j\beta^n}< b,\nonumber\\[2mm]
  &\left.\left|\frac{y}{\alpha^j\beta^m}-
    \frac{y'}{\alpha^j\beta^n}\right|\le\Delta\right\}
    ,\label{eqn.Overlaps1}
\end{align}
which holds for all $m$ and $n$ sufficiently large. Suppose without
loss of generality that $m<n$, and let us derive an estimate for the
number of pairs $(y,y')$ which we are counting on the right hand side
above. If $y\in Y^{(m)}$ and $y'\in Y^{(n)}$ satisfy
\begin{equation}\label{eqn.Overlaps1.5}
    |\beta^{n-m}y-y'|\le \alpha^j\beta^n\Delta,
\end{equation}
then, from the definition \eqref{eqn.Y^nDef} of $Y^{(n)}$, we must
also have that
\begin{equation}\label{eqn.Overlaps2}
|\beta^{n-m}y-y'|\ge 1.
\end{equation}
It is worth pointing out that this was the reason for the introduction
of the sets $Y^{(n)}$. Without the lower bound of a fixed positive
constant here, the next part of the argument would not work. From
Eqs. \eqref{eqn.Overlaps1}-\eqref{eqn.Overlaps2} we now derive that,
for all $m$ and $n$ sufficiently large (depending on $a$ and $b$),
\begin{align*}
  \lambda(\mc{A}_m'\cap\mc{A}_n')&\ll_{r,R}\delta\cdot\left((b-a)\,
    \alpha^{j}\beta^m+1\right)\cdot\left(\alpha^j\beta^n\Delta\right)\\[2mm]
  &\ll (b-a)\,\psi(mJ+j)\,\psi(nJ+j)+\alpha^j\beta^n\delta\Delta\\[2mm]
  &\ll (b-a)\,\psi(mJ+j)\,\psi(nJ+j).
\end{align*}
This implies that there is a constant $K>0$ (depending only on $r$ and
$R$) with the property that
\begin{equation}\label{eqn.LimSupMeas}
  \limsup_{N\rar\infty}\left(\sum_{n\le N}
    \lambda (\mc{A}_n')\right)^2\left(\sum_{m,n\le N}
    \lambda (\mc{A}_m'\cap\mc{A}_n')\right)^{-1}\ge (b-a)K.
\end{equation}
Since the sum of the measures of the sets $\mc{A}_n'$ diverges, by
standard arguments from probability theory (see
\cite[Lemma~2.3]{Harm1998}), it follows that the set of $x$ which fall
in infinitely many of the sets $\mc{A}_n'$ has measure greater than or
equal to $(b-a)K$.\vspace*{.1in}

For the final step of the proof, we could appeal directly to
\cite[Proposition 1]{BereDickVela2006}. For completeness, we provide
the following simple argument. Supposing still that the divergence
condition \eqref{eqn.PsiDiv} holds, let $\mc{W}\subseteq\R$ be the set
of $x\in\R$ for which the inequality \eqref{eqn.PsiIneq} is satisfied
by infinitely many $n\in\N$ and $y\in Y$. If it were the case that
$\lambda (\mc{W}^c)>0$ then, by the Lebesgue density theorem
\cite[Theorem 3.21]{Foll1984}, we could find a point of metric density
$x_0$ of the set $\mc{W}^c$. However, this would imply that
\[\lim_{\epsilon\rar 0^+}\frac{\lambda(\mc{W}\cap B(x_0,\epsilon))}{2\epsilon}=0,\]
which contradicts \eqref{eqn.LimSupMeas}.  Therefore, we conclude that
$\lambda(\mc{W}^c)=0$, thereby completing the proof of our main
result.

\vspace{.2in}

{\footnotesize

\noindent
MB: Fakult\"{a}t f\"{u}r Mathematik, Universit\"{a}t Bielefeld,\\
33501 Bielefeld, Germany.\\
mbaake@math.uni-bielefeld.de\\

\noindent
AH: Department of Mathematics, University of Houston,\\
Houston, TX, United States.\\
haynes@math.uh.edu}


\vspace*{.1in}

\begin{thebibliography}{1}

\bibitem{BaakGrim2013}
M.~Baake, U.~Grimm:
\textit{Aperiodic Order. Vol. 1: A Mathematical Invitation},
Cambridge University Press, Cambridge, 2013.
\vspace*{.1in}

\bibitem{BaakGrim2017}
M.~Baake, A.~Haynes,  D.~Lenz,
\textit{Averaging almost periodic functions along exponential
sequences},
in \textit{Aperiodic Order. Vol.~2:
Crystallography and Almost Periodicity},
Baake M and Grimm U (eds.),
Cambridge University Press, Cambridge (2017),
pp.~343--362; \texttt{arXiv:1704.08120}.
\vspace*{.1in}


\bibitem{BereDickVela2006}
V.~Beresnevich, D.~Dickinson, S.~Velani:
\textit{Measure theoretic laws for lim sup sets},
Memoirs Amer. Math. Soc.  \textbf{179}  (2006) no.\ 846.
\vspace*{.1in}

\bibitem{Buge2012}
Y.~Bugeaud:
\textit{Distribution Modulo One and Diophantine Approximation},
Cambridge University Press, Cambridge, 2012.
\vspace*{.1in}

\bibitem{Foll1984}
G.~B.~Folland:
\textit{Real Analysis. Modern Techniques and Their Applications}, 2nd ed.,
Wiley, New York, 1999.
\vspace*{.1in}

\bibitem{Gall1961}
P.~Gallagher:
\textit{Approximation by reduced fractions},
J. Math. Soc. Japan  \textbf{13} (1961),  342--345.
\vspace*{.1in}


\bibitem{Harm1998}
G.~Harman:
\textit{Metric Number Theory},
Oxford University Press, New York, 1998.
\vspace*{.1in}

\bibitem{Walt1982}
P.~Walters:
\textit{An Introduction to Ergodic Theory,}
Springer, New York, 1982.


\end{thebibliography}
\end{document}